\documentclass[12pt]{article}
\usepackage{wallpaper}

\usepackage{tikz}
\usetikzlibrary{tikzmark,calc}
\usepackage{booktabs}

\definecolor{myColor}{RGB}{0,0,200}

\usepackage{graphicx}
\usepackage{subfig}
\usepackage{amsmath}
\usepackage{amsthm} 
\usepackage{amssymb}
\usepackage[font=footnotesize]{caption} 
\usepackage{subfig} 
\usepackage[noadjust]{cite} 
\usepackage{color}
\usepackage{algorithm} 

\theoremstyle{plain}
 \newtheorem{assumption}{Assumption}

\theoremstyle{definition}

\usepackage{algpseudocode} 
\usepackage{enumerate}
\usepackage[hidelinks,colorlinks=false]{hyperref}
\usepackage[left=0.75in, right=0.75in, top=1.25in, bottom=1in]{geometry}
\usepackage{authblk} 

\usepackage{bm}
\usepackage{tikz}
\usetikzlibrary{calc} 
\usetikzlibrary{shapes} 
\usetikzlibrary{chains}
\usetikzlibrary{fit}
\usetikzlibrary{arrows}
\usetikzlibrary{decorations.text} 
\usetikzlibrary{decorations.markings}
\usetikzlibrary{decorations.pathmorphing} 
\usetikzlibrary{shadows}
\usetikzlibrary{patterns}
\usetikzlibrary{matrix}
\usepackage{pgfplots}
\usepackage[europeanresistors]{circuitikz}
\usepackage[outline]{contour} 
\contourlength{1.5pt}
\graphicspath{{figures/}}

\begin{document}

\title{Technical details of distributed localization}
\author{Xu Fang, Xiaolei Li, Lihua Xie}
\date{}
\maketitle

\subsection*{1. Proof of Theorem $1$}

For any three different
nodes $p_i, p_j, p_k$ in $\mathbb{R}^3$, the condition $\theta_i+\theta_j+\theta_k = \pi$ must hold.
The angle constraints can be rewritten as
\begin{align}\label{para1}
& w_{ik}d_{ik}d_{ij}\cos \theta_i + w_{ki}d_{ik}d_{jk}\cos \theta_k =0,  \\
\label{para2}
& w_{ij}d_{ik}d_{ij}\cos \theta_i + w_{ji}d_{ij}d_{jk}\cos \theta_j =0, \\
\label{para3}
& w_{jk}d_{jk}d_{ij}\cos \theta_j + w_{kj}d_{ik}d_{jk}\cos \theta_k =0,
\end{align}
with $w_{ik}^2+w_{ki}^2 \neq 0$, $w_{ij}^2+w_{ji}^2 \neq 0$, and $w_{jk}^2+w_{kj}^2 \neq 0$. 

First, we introduce \textit{Lemma 7} and \textit{Lemma 8} below for proving \text{Theorem 1}.

\textit{Lemma 7.} $p_i, p_j, p_k$ are non-colinear if the  parameters in (1)-(3) satisfy $
w_{ik}w_{ij}w_{jk}w_{ki}w_{ji}w_{kj} = 0$. 

\begin{proof}
When $
w_{ik}w_{ij}w_{jk}w_{ki}w_{ji}w_{kj} = 0$, without loss of generality,
suppose $w_{ik}=0$, 
since $w_{ik}^2+w_{ki}^2 \neq 0$, we have $\theta_k = \frac{\pi}{2}$, $\theta_i + \theta_j= \frac{\pi}{2}$ from \eqref{para1}. Hence, $p_i, p_j, p_k$ are non-colinear. Similarly, we can prove that $p_i, p_j, p_k$ are non-colinear if the parameter $w_{ij}=0$, or $w_{jk}=0$, or $w_{ki}=0$, or $w_{ji}=0$, or $w_{kj}=0$.
\end{proof}

If $
w_{ik}w_{ij}w_{jk}w_{ki}w_{ji}w_{kj} \neq 0$,  \eqref{para1}-\eqref{para3} can be rewritten as
\begin{align}\label{para21}
& \frac{\cos \theta_i}{\cos \theta_k} = - \frac{w_{ki}d_{jk}}{w_{ik}d_{ij}},  \\
\label{para22}
& \frac{\cos \theta_i}{\cos \theta_j} = - \frac{w_{ji}d_{jk}}{w_{ij}d_{ik}}, \\
\label{para23}
&  \frac{\cos \theta_j}{\cos \theta_k} = - \frac{w_{kj}d_{ik}}{w_{jk}d_{ij}}.
\end{align}

From \eqref{para21} and \eqref{para22}, we have
\begin{align}\label{di}
    & d_{ij} = -\frac{\cos \theta_k}{\cos \theta_i}\frac{w_{ki}}{w_{ik}}d_{jk},\\
    \label{dj}
    & d_{ik} = -\frac{\cos \theta_j}{\cos \theta_i}\frac{w_{ji}}{w_{ij}}d_{jk}.
\end{align}

Note that $\cos \theta_i = \frac{d_{ij}^2+d_{ik}^2-d_{jk}^2}{2d_{ij}d_{ik}}$, $\cos \theta_j = \frac{d_{ij}^2+d_{jk}^2-d_{ik}^2}{2d_{ij}d_{jk}}$, $\cos \theta_k = \frac{d_{ik}^2+d_{jk}^2-d_{ij}^2}{2d_{ik}d_{jk}}$. Combining \eqref{di} and \eqref{dj}, it yields
\begin{equation}\label{dd1}
    \frac{d_{ij}^2+d_{ik}^2-d_{jk}^2}{d_{ik}^2+d_{jk}^2-d_{ij}^2} + \frac{d_{ij}^2+d_{ik}^2-d_{jk}^2}{d_{ij}^2+d_{jk}^2-d_{ik}^2} = -(\frac{w_{ki}}{w_{ik}} +\frac{w_{ji}}{w_{ij}}).
\end{equation}

\textit{Lemma 8.} When the parameters in (1)-(3) satisfy $
w_{ik}w_{ij}w_{jk}w_{ki}w_{ji}w_{kj} \neq 0$,
$p_i, p_j, p_k $ are  colinear if and only if
\begin{equation}
  \frac{w_{ki}}{w_{ik}} +\frac{w_{ji}}{w_{ij}} =1,  \text{or} \ \frac{w_{ij}}{w_{ji}} + \frac{w_{kj}}{w_{jk}} =1, \text{or} \ \frac{w_{ik}}{w_{ki}} +\frac{w_{jk}}{w_{kj}} =1.
\end{equation}
\begin{proof}
(Necessity) If $p_i, p_j, p_k $ are colinear, there are three cases: $\text{(\romannumeral1)}$ $\theta_i=\pi$, $\theta_j, \theta_k =0$; $\text{(\romannumeral2)}$ $\theta_j=\pi$, $\theta_i, \theta_k =0$; $\text{(\romannumeral3)}$ $\theta_k=\pi$, $\theta_i, \theta_j =0$.
For the case $\text{(\romannumeral1)}$ that $\theta_i=\pi$, $\theta_j, \theta_k =0$, we have $d_{ij}+d_{ik}= d_{jk}$. Substituting $d_{ij}+d_{ik}= d_{jk}$ into \eqref{dd1}, we get
\begin{equation}
   \frac{w_{ki}}{w_{ik}} +\frac{w_{ji}}{w_{ij}}=1. 
\end{equation}

Similarly, the conditions can be derived for the other two cases $\text{(\romannumeral2)}$-$\text{(\romannumeral3)}$.


(Sufficiency) If $\frac{w_{ki}}{w_{ik}} +\frac{w_{ji}}{w_{ij}}=1$, \eqref{dd1} becomes
\begin{equation}\label{dd2}
  \frac{d_{ij}^2+d_{ik}^2-d_{jk}^2}{d_{ik}^2+d_{jk}^2-d_{ij}^2} + \frac{d_{ij}^2+d_{ik}^2-d_{jk}^2}{d_{ij}^2+d_{jk}^2-d_{ik}^2} = -1.   
\end{equation}

Then, \eqref{dd2} can be rewritten as
\begin{equation}\label{dd3}
    (d_{ij}^2+d_{ik}^2-d_{jk}^2)^2 = 4d_{ik}^2d_{ij}^2.
\end{equation}

Since $\cos \theta_i = \frac{d_{ij}^2+d_{ik}^2-d_{jk}^2}{2d_{ij}d_{ik}}$, \eqref{dd3} becomes
\begin{equation}
    4d_{ij}^2d_{ik}^2\cos^2 \theta_i =  4d_{ij}^2d_{ik}^2, \rightarrow \cos^2 \theta_i = 1.
\end{equation}

Hence, $\theta_i = 0$ or $\pi$, i.e.,  $p_i, p_j, p_k $ must be colinear. Similarly, we can prove that $p_i, p_j, p_k $ must be colinear for the other two cases $\frac{w_{ij}}{w_{ji}} + \frac{w_{kj}}{w_{jk}} =1 \ \text{and} \ \frac{w_{ik}}{w_{ki}} +\frac{w_{jk}}{w_{kj}} =1$.
\end{proof}


Next, we will prove that the angles $\theta_i, \theta_j, \theta_k \in [0, \pi] $ are determined uniquely by the parameters $w_{ik}, w_{ki}, w_{ij}, w_{ji}, $ $ w_{jk}, w_{kj}$ in (1)-(3). From \textit{Lemma 7} and \textit{Lemma 8}, we can know that
there are only three cases for (1)-(3):

\begin{enumerate}[(i)]
\item $
w_{ik}w_{ij}w_{jk}w_{ki}w_{ji}w_{kj} = 0$;
\item $
w_{ik}w_{ij}w_{jk}w_{ki}w_{ji}w_{kj} \neq 0$, and $\frac{w_{ki}}{w_{ik}} +\frac{w_{ji}}{w_{ij}} =1$,  \text{or} \ $\frac{w_{ij}}{w_{ji}} + \frac{w_{kj}}{w_{jk}} =1$, \text{or} \ $\frac{w_{ik}}{w_{ki}} +\frac{w_{jk}}{w_{kj}} =1$;
\item $
w_{ik}w_{ij}w_{jk}w_{ki}w_{ji}w_{kj} \neq 0$, and $ \frac{w_{ki}}{w_{ik}} +\frac{w_{ji}}{w_{ij}},  \frac{w_{ij}}{w_{ji}} + \frac{w_{kj}}{w_{jk}}, \frac{w_{ik}}{w_{ki}} +\frac{w_{jk}}{w_{kj}} \neq 1$.
\end{enumerate}

The above three cases $(\text{\romannumeral1})\!-\!(\text{\romannumeral3})$ are analyzed below.

For the case $(\text{\romannumeral1})$,
from \textit{Lemma  7}, we can know that $p_i, p_j, p_k$ are non-colinear and form a triangle $\bigtriangleup_{ijk}(p)$.
Without loss of generality,
suppose $w_{ik}=0$, 
since $w_{ik}^2+w_{ki}^2 \neq 0$, we have $\theta_k = \frac{\pi}{2}$, $\theta_i + \theta_j= \frac{\pi}{2}$ from \eqref{para1}. 
Since $\theta_i + \theta_j = \frac{\pi}{2}$, we have
$w_{ij}\cdot w_{ji} < 0$ from \eqref{para22}.  According to the sine rule,
$\frac{d_{jk}}{d_{ik}}= \frac{\sin \theta_i}{\sin \theta_j}$. Then, \eqref{para22} becomes
\begin{equation}\label{tan}
        \frac{\tan \theta_i}{\tan \theta_j} = - \frac{w_{ij}}{w_{ji}}.
\end{equation}

Since $\theta_j= \frac{\pi}{2}-\theta_i$, from \eqref{tan}, we have
\begin{equation}
    \tan \theta_i = \sqrt{- \frac{w_{ij}}{w_{ji}}}, \rightarrow \theta_i = \arctan \sqrt{- \frac{w_{ij}}{w_{ji}}}.
\end{equation}

Similarly, we can prove that $\theta_i, \theta_j, \theta_k$ can be determined uniquely if the parameter $w_{ij}, w_{jk}, w_{ki}, w_{ji},$ or $w_{kj}$ equals $0$.

For the case $(\text{\romannumeral2})$, from \textit{Lemma 8}, we can know that $p_i, p_j, p_k$ are colinear. Two of $\theta_i, \theta_j, \theta_k$ must be $0$. If $w_{kj}w_{jk}<0$, i.e., $\frac{\cos \theta_j}{\cos \theta_k}>0$ from \eqref{para23}, 
we have
$\theta_i=\pi$, $\theta_j, \theta_k =0$. Similarly, we have $\theta_j=\pi$, $\theta_i, \theta_k =0$ if $w_{ki}w_{ik}<0$, and
$\theta_k=\pi$, $\theta_i, \theta_j =0$ if $w_{ji}w_{ij}<0$.

For the case $(\text{\romannumeral3})$, from \textit{Lemma 8}, we can know that  $p_i, p_j, p_k$ are non-colinear and form a triangle $\bigtriangleup_{ijk}(p)$. For this triangle $\bigtriangleup_{ijk}(p)$,
at most one of $\theta_i, \theta_j, \theta_k$ is an obtuse angle. Hence,
there are only four possible cases: $(\text{a})$  $w_{ki}w_{ik}, w_{ji}w_{ij}, w_{kj}w_{jk} <0$; $(\text{b})$ $w_{ki}w_{ik}, w_{ji}w_{ij}>0, w_{kj}w_{jk} <0$; $(\text{c})$ $w_{ki}w_{ik}, w_{kj}w_{jk}>0, w_{ji}w_{ij}<0$; $(\text{d})$ $w_{ji}w_{ij}, w_{kj}w_{jk} >0, w_{ki}w_{ik}<0$. 
For the case $(\text{a})$, we have $\theta_i, \theta_j, \theta_k < \frac{\pi}{2}$. From \eqref{para21} and \eqref{para22}, we have
\begin{equation}\label{trans}
    \tan \theta_k = -\frac{w_{ki} }{w_{ik}} \tan \theta_i, \ \ \tan \theta_j = -\frac{w_{ji} }{w_{ij}} \tan \theta_i. 
\end{equation}

Note that $\tan \theta_i = \tan (\pi- \theta_j - \theta_k)= \frac{\tan \theta_j +\tan \theta_k}{\tan \theta_j\tan \theta_k -1}$. Based on \eqref{trans}, 
we have
\begin{equation}
\tan \theta_i = \sqrt{\frac{1-\frac{w_{ki}}{w_{ik}}- \frac{w_{ji}}{w_{ij}}}{\frac{w_{ki}w_{ji}}{w_{ik}w_{ij}}}}.
\end{equation}

Then, we can obtain the angle $\theta_i$ by
\begin{equation}
    \theta_i = \arctan \sqrt{\frac{1-\frac{w_{ki}}{w_{ik}}- \frac{w_{ji}}{w_{ij}}}{\frac{w_{ki}w_{ji}}{w_{ik}w_{ij}}}}.
\end{equation}

Similarly, the angles $\theta_j$ and $\theta_k$ can also be obtained. Using this way, we can prove that $\theta_i, \theta_j, \theta_k$ can be determined uniquely by the parameters $w_{ik}, w_{ki}, w_{ij}, w_{ji}, w_{jk}, w_{kj}$ for the cases $(\text{b})$-$(\text{d})$.

\begin{figure}[t]
\centering
\includegraphics[width=0.4\linewidth]{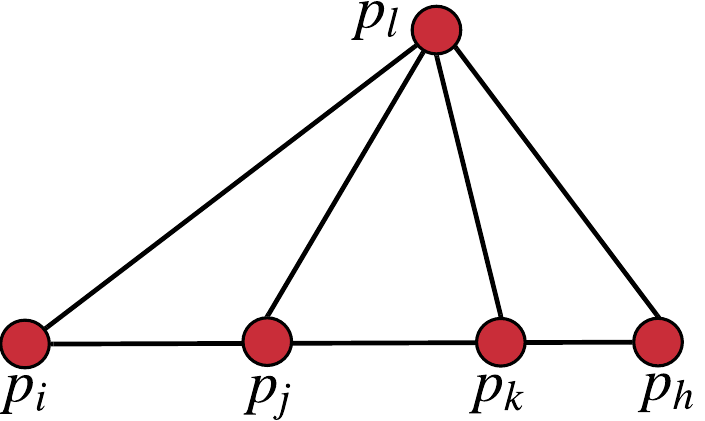}
\caption{3-D local-relative-bearing-based network. }
\label{moti1}
\end{figure}

\subsection*{2.  Proof of \textit{Lemma} $2$}

\begin{proof}
Since $ \mu_{ij}e_{ij}\!+\! \mu_{ik}e_{ik}\!+\! \mu_{ih}e_{ih}\!+\! \mu_{il}e_{il}=\mathbf{0}$ and $w_{ik}e_{ik}^Te_{ij}+w_{ki}e_{ki}^Te_{kj}=0$,
for the scaling space $S_s$,  it is straightforward that $\eta_d^Tp =\mathbf{0}$ and $\eta_r^Tp =0$. For the translation space $S_t$, we have $\eta_d^T( \mathbf{1}_n\otimes {I}_3)=\mathbf{0}$ and $\eta_r^T( \mathbf{1}_n\otimes {I}_3)=\mathbf{0}$.  For the rotation space $S_r= \{(I_n \otimes A)p, A+A^T =\mathbf{0}, A \in \mathbb{R}^{3 \times 3} \}$, it follows that $\eta_d^T (I_n \otimes A) p = A (\mu_{ij}e_{ij}\!+\! \mu_{ik}e_{ik}\!+\! \mu_{ih}e_{ih}\!+\! \mu_{il}e_{il})=\mathbf{0}$ and 
\begin{equation}
\begin{array}{ll}
\eta_r^T (I_n \otimes A) p\\
= w_{ik}p_i^T(A+A^T)p_i \!+\! (w_{ki}\!-\!w_{ik})p_j^T(A\!+\!A^T)p_i 
\\
-(w_{ik}\!+\!w_{ki})p_k^T(A\!+\!A^T)p_i + (w_{ik}\!-\!w_{ki})p_k^T(A\!+\!A^T)p_j \\
+w_{ki}p_k^T(A\!+\!A^T)p_k =0.
\end{array}
\end{equation}

Then, the conclusion follows.
\end{proof}

\subsection*{3.  Local-relative-bearing-based Displacement Constraint}

$g_{ij} = \frac{e_{ij}}{d_{ij}} \in \mathbb{R}^3$ is the relative bearing of $p_j$ with respect to $p_i$ in $\Sigma_g$. For the node $i$ and its neighbors $j, k, h, l$ in $\mathbb{R}^3$,
the matrix $g_i=(g_{ij}, g_{ik}, g_{ih}, g_{il}) \in \mathbb{R}^{3 \times 4}$ is a wide matrix. From the matrix theory, there must be a non-zero vector $ \bar \mu_i=(\bar \mu_{ij}, \bar \mu_{ik}, \bar \mu_{ih}, \bar \mu_{il})^T \in \mathbb{R}^4$ such that $g_i\bar \mu_i= \mathbf{0}$, i.e.,

\begin{equation}\label{root1}
     \bar \mu_{ij}g_{ij}+\bar \mu_{ik}g_{ik}+\bar \mu_{ih}g_{ih} + \bar \mu_{il}g_{il}= \mathbf{0},
\end{equation}
where $\bar \mu_{ij}^2+\bar \mu_{ik}^2+\bar \mu_{ih}^2+\bar \mu_{il}^2 \neq 0$.

The equation $g_i\bar \mu_i= \mathbf{0}$ is a bearing constraint, based on which a displacement constraint  can be obtained shown as following. The non-zero vector $(\bar \mu_{ij}, \bar \mu_{ik}, \bar \mu_{ih}, \bar \mu_{il})^T$ can be calculated with local relative bearing measurements $g_{ij}^{i}, g_{ik}^{i}, g_{ih}^{i}, g_{il}^{i}$ by solving the following equation
\begin{equation}\label{wmi1}
\left[ \!
\begin{array}{c c c c}
g_{ij}^{i} & g_{ik}^{i}  & g_{ih}^{i} &  g_{il}^{i} \\
\end{array}
\right]  \left[ \!
	\begin{array}{c}
	\bar \mu_{ij} \\
	\bar \mu_{ik} \\
	\bar \mu_{ih} \\
	\bar \mu_{il}
	\end{array}
	\right] = \mathbf{0}.
\end{equation}

Note that \eqref{root1} can be rewritten as
\begin{equation}\label{loca1}
\bar \mu_{ij}\frac{e_{ij}}{d_{ij}}+\bar \mu_{ik}\frac{e_{ik}}{d_{ik}}+\bar \mu_{ih}\frac{e_{ih}}{d_{ih}} + \bar \mu_{il}\frac{e_{il}}{d_{il}} = \mathbf{0}.
\end{equation}

\begin{assumption}\label{ad3}
No two nodes are collocated in $\mathbb{R}^3$. Each anchor node has at least two neighboring anchor nodes, and
each free node has at least four neighboring nodes. The free node and its neighbors are non-colinear.
\end{assumption}

Under \text{Assumption \ref{ad3}}, without loss of generality,  suppose node $l$ is not colinear with nodes $i,j,k,h$  shown in the above Fig. \ref{moti1}.
The angles among the nodes $p_i, p_j, p_k, p_h, p_l$ are denoted by $\xi_{ilj} \!=\! \angle p_ip_lp_j, \xi_{ilk} \!=\! \angle p_ip_lp_k, \xi_{ilh} \!=\! \angle p_ip_lp_h, \xi_{ijl} \!=\! \angle p_ip_jp_l, \xi_{ikl} \!=\! \angle p_ip_kp_l, \xi_{ihl} \!=\! \angle p_ip_hp_l$. Note that these angles can be obtained by only using the local relative bearing measurements. For example, $\xi_{ilj} = g_{li}^Tg_{lj}= {g^{l}_{li}}^TQ_l^TQ_lg_{lj}^{l} ={g_{li}^{l}}^Tg_{lj}^l$. 
According to the sine rule,
$\frac{d_{il}}{d_{ij}}= \frac{\sin \xi_{ijl}}{\sin \xi_{ilj}}, \frac{d_{il}}{d_{ik}}= \frac{\sin \xi_{ikl}}{\sin \xi_{ilk}}, \frac{d_{ih}}{d_{il}}= \frac{\sin \xi_{ilh}}{\sin \xi_{ihl}}$. Then, based on \eqref{loca1},
we can obtain
a displacement constraint by only using the local relative bearing measurements shown as
\begin{equation}\label{bead}
     \mu_{ij}e_{ij}+ \mu_{ik}e_{ik}+  \mu_{ih}e_{ih} +  \mu_{il}e_{il}= \mathbf{0},
\end{equation}
where
\begin{equation}
\begin{array}{ll}
     &  \mu_{ij} = \bar \mu_{ij}\frac{\sin \xi_{ijl}}{\sin \xi_{ilj}}, \ \  \mu_{ik} = \bar \mu_{ik}\frac{\sin \xi_{ikl}}{\sin \xi_{ilk}}, \\
     & \mu_{ih} = \bar \mu_{ih}\frac{\sin \xi_{ilh}}{\sin \xi_{ihl}}, \ \  \mu_{il} = \bar \mu_{il}.
\end{array}
\end{equation}

In a local-relative-bearing-based network in $\mathbb{R}^3$ under \text{Assumption \ref{ad3}},
let 
$\mathcal{X}_{\mathcal{G}}= \{ ( i, j, k, h, l) \in \mathcal{V}^{5} : (i,j), (i,k), $ $ (i,h),  (i,l),  (j,k), (j,h),   (j,l)  \in \mathcal{E},  j \!<\! k \!<\! h \!<\! l\}$. Each element of $\mathcal{X}_{\mathcal{G}}$ can be used to construct a local-relative-bearing-based displacement constraint.

\subsection*{4.  Distance-based Displacement Constraint}

Since the displacement constraints
are invariant to translations and rotations, a congruent network of the subnetwork consisting of the node and its neighbors
has the displacement constraint. 
Each displacement constraint can be regarded as a subnetwork, and
multi-dimensional scaling can be used to obtain 
displacement constraint shown in the following Algorithm \ref{disa} \cite{han2017barycentric}.

\subsection*{5. Ratio-of-distance-based Displacement Constraint}

For the free node $i$ and its neighbors $j,k,h,l$, under Assumption $1$, we can obtain the ratio-of-distance matrix $M_r$ \eqref{ratio} by the ratio-of-distance measurements.

\begin{equation}\label{ratio}
M_r = \frac{1}{d_{ij}^2}\left[ \!
\begin{array}{c c c c c}
0 & d_{ij}^2 & d_{ik}^2  & d_{ih}^2 &  d_{il}^2 \\
d_{ji}^2 & 0 & d_{jk}^2 & d_{jh}^2 & d_{jl}^2 \\
d_{ki}^2 & d_{kj}^2 & 0 & d_{kh}^2 & d_{kl}^2 \\
d_{hi}^2 & d_{hj}^2 & d_{hk}^2 & 0 & d_{hl}^2 \\
d_{li}^2 & d_{lj}^2 & d_{lk}^2 & d_{lh}^2 & 0
\end{array}
\right].   
\end{equation}

Note that the displacement constraints
are not only invariant to translations and rotations, but also to scalings.  Hence,  a network with ratio-of-distance measurements  $\frac{1}{d_{ij}}\{d_{ij},  \cdots, d_{hl}, \cdots \}$ has the same displacement constraints as the network with  distance measurements $\{d_{ij},  \cdots, d_{hl}, \cdots \}$, that is, the displacement constraint $\mu_{ij}e_{ij}+ \mu_{ik}e_{ik}+  \mu_{ih}e_{ih} +  \mu_{il}e_{il}= \mathbf{0}$ can also be obtained by Algorithm $1$, where the distance matrix  $M$ \eqref{dim} is replaced by the
the ratio-of-distance matrix $M_r$ \eqref{ratio}. 

\begin{algorithm}
\caption{Distance-based displacement constraint}
\label{disa}
\begin{algorithmic}[1]
\State Available information: Distance measurements among the nodes $p_i,p_j,p_k,p_h,p_l$. Denote  $(\mathcal{\bar G}, \bar p)$ as a subnetwork with $\bar p=(p_i^T,p_j^T,p_k^T,p_h^T,p_l^T)^T$.
\State Constructing a distance matrix $M \in \mathbb{R}^{5 \times 5}$ shown as
\begin{equation}\label{dim}
    M = \left[ \!
\begin{array}{c c c c c}
0 & d_{ij}^2 & d_{ik}^2  & d_{ih}^2 &  d_{il}^2 \\
d_{ji}^2 & 0 & d_{jk}^2 & d_{jh}^2 & d_{jl}^2 \\
d_{ki}^2 & d_{kj}^2 & 0 & d_{kh}^2 & d_{kl}^2 \\
d_{hi}^2 & d_{hj}^2 & d_{hk}^2 & 0 & d_{hl}^2 \\
d_{li}^2 & d_{lj}^2 & d_{lk}^2 & d_{lh}^2 & 0
\end{array}
\right]. 
\end{equation}
\State Computer the centering matrix $J=I-\frac{1}{5}\mathbf{1}_5\mathbf{1}_5^T$;
\State Compute the matrix $X=-\frac{1}{2}JMJ$;
\State Perform singular value decomposition on $X$ as
\begin{equation}
    X = V \Lambda V^T,
\end{equation}
where $V=(v_1,v_2,v_3,v_4,v_5) \in \mathbb{R}^{5 \times 5}$ is a unitary matrix, and $\Lambda = \text{diag}(\lambda_1,\lambda_2,\lambda_3, \lambda_4, \lambda_5)$ is a diagonal matrix whose diagonal elements $\lambda_1 \ge \lambda_2 \ge \lambda_3 \ge \lambda_4 \ge \lambda_5$ are singular values. Since $\text{Rank}(X) \le 3$, we have $\lambda_4 = \lambda_5=0$. Denote by $V_*=(v_1,v_2,v_3)$ and $\Lambda_*=\text{diag}(\lambda_1,\lambda_2,\lambda_3)$;
\State 
Obtaining a congruent network $(\mathcal{\bar G}, \bar q) \cong (\mathcal{\bar G}, \bar p)$ with $\bar q=(q_i^T,q_j^T,q_k^T,q_h^T,q_l^T)^T $, where $(q_i, q_j, q_k, $ $ q_h, q_l ) = \Lambda_* ^{\frac{1}{2}}V_*^T$;
\State Based on the congruent network $\bar q=(q_i^T,q_j^T,q_k^T,q_h^T,q_l^T)^T$ of the subnetwork $\bar p=(p_i^T,p_j^T,p_k^T,p_h^T,p_l^T)^T$, the parameters $\mu_{ij}, \mu_{ik}, \mu_{ih}, \mu_{il}$ in $\mu_{ij}e_{ij}+\mu_{ik}e_{ik}+\mu_{ih}e_{ih} + \mu_{il}e_{il}= \mathbf{0}$ can be obtained by solving the following matrix equation
\begin{equation}\label{wmi}
\left[ \!
\begin{array}{c c c c}
q_j-q_i & q_k-q_i  & q_h-q_i &  q_l-q_i \\
\end{array}
\right]  \left[ \!
	\begin{array}{c}
	\mu_{ij} \\
	\mu_{ik} \\
	\mu_{ih} \\
	\mu_{il}
	\end{array}
	\right] = \mathbf{0}.
\end{equation}
\end{algorithmic}
\end{algorithm}

\subsection*{6. Angle-based Displacement Constraint}

For a triangle $\bigtriangleup_{ijk}(p)$, according to the sine rule,
the ratios of distance 
can be calculated by the angle measurements $\theta_i, \theta_j, \theta_k$ shown as
\begin{equation}\label{sin}
\frac{d_{ij}}{d_{ik}}=\frac{\sin \theta_k}{\sin \theta_j}, \frac{d_{ij}}{d_{jk}}=\frac{\sin \theta_k}{\sin \theta_i}.
\end{equation}

Under \text{Assumption \ref{ad3}},  the ratios of distance of all the edges among the nodes $i, j, k, h, l$ can be calculated by the angle measurements through the sine rule \eqref{sin}, i.e., the ratio-of-distance matrix $M_r$ \eqref{ratio} is available. Then, the displacement constraint $\mu_{ij}e_{ij}+ \mu_{ik}e_{ik}+  \mu_{ih}e_{ih} +  \mu_{il}e_{il}= \mathbf{0}$ can be obtained by Algorithm $1$, where the distance matrix  $M$ \eqref{dim} is replaced by the
the ratio-of-distance matrix $M_r$.

In an angle-based network in $\mathbb{R}^3$ under \text{Assumption \ref{ad3}}, 
let 
$\mathcal{X}_{\mathcal{G}}=\{ ( i, j, k, h, l) \in \mathcal{V}^{5} : (i,j), (i,k), $ $ (i,h),  (i,l),  (j,k), (j,h),  (j,l), (k,h), (k,l), (h,l) \in \mathcal{E}, j \!<\! k \!<\! h \!<\! l\}$. Each element of $\mathcal{X}_{\mathcal{G}}$ can be used to construct an angle-based displacement constraint.

\subsection*{7. Relaxed Assumptions for Constructing local-relative-position-based, Distance-based,  Ratio-of-distance-based, Local-relative-bearing-based, and Angle-based Displacement Constraint in a Coplanar Network}

\begin{assumption}\label{as31}
No two nodes are collocated in $\mathbb{R}^3$. Each anchor node has at least two neighboring anchor nodes, and
each free node has at least three neighboring nodes.
\end{assumption}

\begin{assumption}\label{ad32}
No two nodes are collocated in $\mathbb{R}^3$. Each anchor node has at least two neighboring anchor nodes, and
each free node has at least three neighboring nodes. The free node and its neighbors are non-colinear.
\end{assumption}


\begin{enumerate}

\item 

In a local-relative-position-based coplanar network in $\mathbb{R}^3$ with \text{Assumption \ref{as31}}, 
let 
$\mathcal{X}_{\mathcal{G}}=\{ ( i, j, k, h) \in \mathcal{V}^{4} : (i,j), (i,k), $ $ (i,h)  \in \mathcal{E}, j \!<\! k \!<\! h \}$. Each element of $\mathcal{X}_{\mathcal{G}}$ can be used to construct a local-relative-position-based displacement constraint $\mu_{ij}e_{ij}+\mu_{ik}e_{ik}+\mu_{ih}e_{ih}= \mathbf{0}$.

\item In a distance-based coplanar network in $\mathbb{R}^3$ with \text{Assumption \ref{as31}}, 
let 
$\mathcal{X}_{\mathcal{G}}=\{ ( i, j, k, h) \in \mathcal{V}^{4} : ((i,j), (i,k), $ $ (i,h),  (j,k),  (j,h),   (k,h) \in \mathcal{E}, j \!<\! k \!<\! h\}$. Each element of $\mathcal{X}_{\mathcal{G}}$ can be used to construct a distance-based displacement constraint $\mu_{ij}e_{ij}+\mu_{ik}e_{ik}+\mu_{ih}e_{ih}= \mathbf{0}$. 
\item In a ratio-of-distance-based coplanar network in $\mathbb{R}^3$ with \text{Assumption \ref{as31}},
let 
$\mathcal{X}_{\mathcal{G}}=\{ ( i, j, k, h) \in \mathcal{V}^{4} : ((i,j), (i,k), $ $ (i,h),  (j,k),  (j,h),   (k,h) \in \mathcal{E}, j \!<\! k \!<\! h\}$. Each element of $\mathcal{X}_{\mathcal{G}}$ can be used to construct a ratio-of-distance-based displacement constraint $\mu_{ij}e_{ij}+\mu_{ik}e_{ik}+\mu_{ih}e_{ih}= \mathbf{0}$.
\item In a local-relative-bearing-based coplanar network in $\mathbb{R}^3$ with \text{Assumption \ref{ad32}}, 
let 
$\mathcal{X}_{\mathcal{G}}=\{ ( i, j, k, h)  $ $ \in \mathcal{V}^{4} : (i,j), (i,k), (i,h),  (j,k), (j,h) \in \mathcal{E}, j \!<\! k \!<\! h \!<\! l\}$. Each element of $\mathcal{X}_{\mathcal{G}}$ can be used to construct a local-relative-bearing-based displacement constraint $\mu_{ij}e_{ij}+\mu_{ik}e_{ik}+\mu_{ih}e_{ih}= \mathbf{0}$.
\item In an angle-based coplanar network in $\mathbb{R}^3$ with \text{Assumption \ref{ad32}}, 
let 
$\mathcal{X}_{\mathcal{G}}=\{ ( i, j, k, h) \in \mathcal{V}^{4} : ((i,j), (i,k), $ $ (i,h),  (j,k),  (j,h),   (k,h) \in \mathcal{E}, j \!<\! k \!<\! h\}$. Each element of $\mathcal{X}_{\mathcal{G}}$ can be used to construct an angle-based displacement constraint $\mu_{ij}e_{ij}+\mu_{ik}e_{ik}+\mu_{ih}e_{ih}= \mathbf{0}$.
\end{enumerate}

\bibliographystyle{IEEEtran}
\bibliography{papers}

\end{document}